\documentstyle[11pt,twoside]{amsart}
\title{Moduli spaces $M_{g,n}(W)$ for surfaces}
\author{Valery Alexeev}
\date{October 3, 1994}

\newcommand{\cX}{{\cal X}}
\newcommand{\cY}{{\cal Y}}
\newcommand{\cR}{{\cal R}}
\newcommand{\cS}{{\cal S}}
\newcommand{\cT}{{\cal T}}
\newcommand{\cU}{{\cal U}}
\newcommand{\cV}{{\cal V}}
\newcommand{\cB}{{\cal B}}
\newcommand{\cC}{{\cal C}}
\newcommand{\cE}{{\cal E}}
\newcommand{\cO}{{\cal O}}
\newcommand{\cL}{{\cal L}}
\newcommand{\cM}{{\cal M}}
\newcommand{\cA}{{\cal A}}
\newcommand{\cF}{{\cal F}}
\newcommand{\cH}{{\cal H}}

\newcommand{\cGamma}{{\cal\Gamma}}

\newcommand{\cMC}{{{\cM\cC}} }
\newcommand{\bMC}{{{\bold M \bold C}} }
\newcommand{\cMCS}{{$\cM\cC(\cS)$} }

\newcommand{\codim}{\operatorname{codim} }

\newcommand{\Aut}{\operatorname{Aut} }

\newcommand{\Sym}{\operatorname{Sym} }

\newtheorem{thm}{Theorem}[section]
\newtheorem{question}[thm]{Question}
\newtheorem{lem}[thm]{Lemma}

\theoremstyle{definition}
\newtheorem{defn}[thm]{Definition}

\newtheorem{say}[thm]{}
\newtheorem{exmp}[thm]{Example}

\newtheorem{rem}[thm]{Remark}           
\newtheorem{ack}{Acknowledgments}       
\newtheorem{notation}{Notation}         
\newtheorem{notationnum}[thm]{Notation}

\theoremstyle{remark}

\newtheorem{approach}{Approach}
\newtheorem{principle}{Principle}

\begin{document}
\bibliographystyle{amsplain}
\maketitle
\tableofcontents

\setcounter{section}{-1}
\section{Introduction}
\label{sec:introduction}

\begin{say}
  \label{say:moduli spaces for curves}
  We recall the following coarse moduli spaces in the case of curves:
  \begin{enumerate}
  \item $M_{g}$, parameterizing nonsingular curves of genus $g\ge2$
    and its compactification $\overline{M_{g}}$, parameterizing
    Mumford-Deligne moduli-stable curves, see Mumford \cite{Mumford77},
  \item spaces $M_{g,n}$, $2g-2+n>0$, for stable
    $n$-pointed curves, see Knudsen \cite{Knudsen83},
  \item a moduli space $M_{g,n}(W)$ of stable maps from reduced curves
    to a variety $W$, see Kontsevich \cite{Kontsevich94}.
  \end{enumerate}
  It is well known that $\overline{M_{g}}$ and $M_{g,n}$ are
  projective, $M_{g}$ is quasi-projective.
\end{say}

\begin{say}
  For surfaces, Gieseker \cite{Gieseker77} established the existence
  of a quasi-projective scheme parameterizing surfaces with at worst
  Du Val singularities, ample canonical class $K$ and fixed $K^{2}$,
  this is a straightforward analog of $M_{g}$ and we will denote it
  by $M_{K^{2}}$. A geometrically meaningful compactification of this
  space, $\overline{M^{sm}_{K^{2}}}$, was constructed by Koll\'ar and
  Shepherd-Barron in \cite{KollarShepherdBarron88} as a separated
  algebraic space.  It is a moduli space of smoothable stable (not in
  the G.I.T. sense) surfaces of general type.  In \cite{Kollar90}
  Koll\'ar has shown that if the class of smoothable stable surfaces
  with a fixed $K^2$ is bounded then $\overline{M^{sm}_{K^{2}}}$ is in
  fact a projective scheme (Corollary 5.6). Finally, the boundedness
  was proved by the author in \cite{Alexeev94b}.
\end{say}

\begin{say}
  The main purpose of this paper is to construct analogs of $M_{g,n}$
  and $M_{g,n}(W)$ in the case of surfaces, and to prove their
  projectiveness. After this is done, we touch on a connection between
  our moduli spaces and the standard moduli spaces of K3 and Abelian
  surfaces.
\end{say}

\begin{say}
  An idea of ``$M_{g,n}$ for surfaces'' occured to me when I mentioned
  that my boundedness theorem 9.2 \cite{Alexeev94b} is strictly
  stronger than what was used for $\overline{M^{sm}_{K^{2}}}$. Then,
  looking at the definition of $M_{g,n}(W)$ in \cite{Kontsevich94} I
  realized that this is simply a relative version of the same scheme,
  and can be done for surfaces too.
\end{say}

\begin{say}
  The basic construction of a moduli space as an algebraic space used
  here is the same as in \cite{Mumford82}, \cite{Kollar85},
  \cite{Viehweg94} and elsewhere. For the hardest question involved,
  proof of local closedness, we refer to a result of Koll\'ar
  \cite{Kollar94}.

  For proving that our moduli spaces are projective schemes, rather
  than mere algebraic spaces, we use Koll\'ar's Ampleness Lemma 3.9
  \cite{Kollar90}, which can be applied in a straightforward way to a
  variety of complete moduli problems.
\end{say}

\begin{say}
  Kontsevich and Manin \cite{KontsevichManin94} use the moduli spaces
  $M_{g,n}(W)$ to define Gromov-Witten classes of varieties in the
  ``quantum cohomology'' theory. Hence one distant application of
  ``$M_{g,n}(W)$ for surfaces'' might be ``higher'' GW-classes of
  schemes.
\end{say}

\begin{notation}
  All schemes are assumed to be at least Noetherian and defined over a
  fixed algebraically closed field $k$ of characteristic zero.
  Obstacles to extending the results to positive characteristic are
  discussed in the last section.  In most cases, we prefer the
  additive notation to the multiplicative one, for divisors and line
  bundles. All moduli spaces in this paper are coarse.
\end{notation}

\begin{ack}
  It is a pleasure to acknowledge useful discussions that I had with
  F.Campana, Y.Kawamata, J.Koll\'ar, Yu.Manin, E.Sernesi and
  V.V.Shokurov while working on this paper.
\end{ack}

\section{Overview}
\label{sec:overview}

\begin{say}
  One possible approach to solving an algebro-geometric moduli
  problem goes through the following steps:
  \begin{enumerate}
  \item defining the objects that we are trying to parameterize,
  \item giving the right definition for a moduli functor,
  \item establishing properties of this functor,
  \item constructing a moduli space in some fashion,
  \item proving finer facts about this space.
  \end{enumerate}
  In our treatment, we will follow two guiding principles well
  understood in the field:
  \begin{principle}
    Most moduli spaces exist in the category of algebraic spaces.
  \end{principle}
  \begin{principle}
    Most complete and separated moduli spaces are projective.
  \end{principle}
\end{say}

\begin{say}
  Moduli spaces of nonsingular curves $M_g$ and of Deligne-Mumford
  stable curves $\overline{M_{g}}$ of genus $g$ provide a textbook
  illustration of how this works in practice. Since nonsingular curves
  can degenerate into singular ones in an uncorrectable way, $M_g$ is
  not complete. There are many different ways to compactify it but the
  one we are interested in here is adding more curves and trying to
  solve an enlarged moduli problem. It turns out that the curves one
  has to add are Deligne-Mumford moduli-stable curves which are
  defined as connected and complete reduced curves with ordinary nodes
  only such that every smooth rational irreducible component
  intersects others in at least 3 points and every irreducible
  component of arithmetical genus one intersects the rest in at least
  1 point. The latter condition has two equivalent meanings:
  \begin{enumerate}
    \label{enu:properties of deligne-mumford curves}
  \item the automorphism group $\Aut(X)$ is finite (and this property
    is a must for the Geometric Invariant Theory),
  \item the dualizing sheaf $\omega_X$ is ample.
  \end{enumerate}
  To arrive at this answer, one can look at the way the good
  limits are obtained. One considers a family $\cX$ of curves over
  a marked curve, or the specter of a DVR ring, $(\cS,0)$ with a
  nonsingular general fibre and a degenerate special fibre. Then by
  the Semistable Reduction Theorem, after making a finite base change
  $\cS'\to\cS$ and resolving the singularities of
  $\cX'$, the central fibre will be a reduced curve with
  simple nodes. Following (1) above one should contract all the
  rational curves $E$ in the central fiber that intersect the rest
  only at 1 or 2 points. These have self-intersection numbers $E^2=-1$
  and $E^2=-2$ respectfully. Contracting $(-1)$-curves leaves the
  ambient space, which is a surface, nonsingular. Contracting
  $(-2)$-curves introduces very simple surface singularities, called
  Du Val or rational double. The central fiber has nodes only.
\end{say}

\begin{say}
  One can recognize that the above is a field of the Minimal Model
  Program. In fact, we have just constructed the canonical model, in
  dimension 2, of $\cX'$ over $\cS'$. So, to generalize
  $\overline{M_{g}}$ to the surfaces of general type we have to
  apply the Minimal Model Program in dimension 3. This was done by
  Koll\'ar and Shepherd-Barron in \cite{KollarShepherdBarron88}. By
  that time, the end of 1980-s, all the necessary for this
  construction tools from MMP in dimension 3 were already available.
  The new objects that one has to add are defined as connected reduced
  surfaces with semi-log canonical singularities and ample tensor
  power of the dualizing sheaf $(\omega_X^N)^{**}$, where $^{**}$
  means taking the self-dual. Using the additive notation, we say that
  a $\Bbb Q$-Cartier divisor $K_X$ is ample.
\end{say}

\begin{say}
  At the present time, the {\em log} Minimal Model Program in
  dimension 3 is in a pretty good shape. Let us see what kind of
  statements we can get using its principles.  Keeping in line with
  what we did before, we now consider pairs $(X,B)$ of surfaces $X$
  and divisors $B=\sum_{j=1}^{n} B_j$ with ample $K_X+B$. A
  construction very similar to the one above, with Semistable
  Reduction Theorem and, this time {\em log} canonical model, shows
  that we again get a complete moduli functor. What about
  singularities of the pair $(X,B)$? Why, they ought to be semi-log
  canonical, of course!

  What is the analog of this in dimension one? That is easy to answer
  and we get something very familiar. The divisor
  $B=\sum_{j=1}^{n} B_j$ becomes a set of marked points. Semi-log
  canonical means that the curve has nodes only, and that marked
  points are distinct and lie in the nonsingular part. These are
  exactly the $n$-marked semistable curves of Knudsen
  \cite{Knudsen83}.
\end{say}

\begin{say}
  Another possible generalization would be looking not at absolute
  curves (or surfaces) $X$ (or pairs $(X,B)$) but doing it in the
  relative setting. In other words, let us consider maps $X\to W$ to a
  fixed projective scheme $W$ with $K_X$ (resp. $K_X+B$) {\em
    relatively\/} ample. The only modification in the above
  construction will be that we have to apply the relative version of
  the (log) Minimal Model Program over $\cS'\times W$ instead of over
  $\cS'$. What we get for curves is the moduli space $M_{g,n}(W)$
  introduced by Kontsevich in \cite{Kontsevich94}.
\end{say}

\begin{say}
  Now that we have outlined the objects we will be dealing with, let
  us return back to the basic example of $\overline{M_g}$. We recall
  two different approaches to constructing it.
\end{say}

\begin{approach}[G.I.T.]
  One first proves that moduli-stable curves are asymptotically
  Hilbert-stable, \cite{Mumford77}. Then the standard G.I.T. machinery
  produces a quasi-projective moduli space. Since it is complete, it is
  actually projective.
\end{approach}

\begin{approach}
  Using a fairly general argument (\cite{Mumford82}, p.172) one proves
  the existence of a moduli space in the category of algebraic spaces.
  To a family of curves $\pi:\cX\to\cS$ one can in a natural way
  associate line bundles on $\cS$ which are defined as
  $\det(\pi_*\omega^k)$. They descend to ($\Bbb Q$-)line bundles
  $\lambda_k$ on $\overline{M_g}$ and one can further show that
  $\lambda_k$ are ample for $k\ge1$.
\end{approach}

\begin{say}
  As mentioned in \cite{ShepherdBarron83} and \cite{Kollar90}, for
  surfaces of general type the first approach fails. By
  \cite{Mumford77} 3.19 in order to be asymptotically Chow- or
  Hilbert-stable a surface has to have singularities of multiplicities
  at most 7. On the other hand, the semi-stable limits described above
  have semi-log canonical singularities and it looks like they must be
  included in any reasonable complete moduli problem. These semi-log
  canonical singularities include all quotient singularities, for
  example, and can have arbitrarily high multiplicities.
\end{say}

\begin{say}
  The second approach is what we will be using here. After
  establishing the existence as an algebraic space, we will use
  Koll\'ar's Ampleness Lemma \cite{Kollar90} to prove that it is
  projective. The Ampleness Lemma is a general scheme that shows
  projectiveness once some good properties of the moduli functor are
  established: local closedness, completeness, separateness,
  semipositiveness, finite reduced automorphism groups, and,
  crucially, boundedness. The projectiveness will be the only
  ``finer'' property of the obtained moduli spaces that we will
  consider.
\end{say}

\section{The objects}
\label{sec:the objects}

\begin{say}
  The main objects into the consideration will be
  {\em stable maps of pairs\/} $g:(X,B)\to W$, where
  \begin{enumerate}
  \item $W\subset\Bbb P$ is a fixed projective scheme,
  \item $X$ is a connected projective surface,
  \item $B=\sum_{j=1}^n B_j$ is a divisor on $X$, $B_j$ are reduced
    but not necessarily irreducible,
  \item the pair $(X,B)$ has semi-log canonical singularities,
  \item the divisor $K_X+B$ is relatively $g$-ample.
  \end{enumerate}
  The precise definitions follow.
\end{say}

\begin{say}
  For a normal variety $X$, $K_{X}$ or simply $K$ will always denote
  the class of linear equivalence of the canonical Weil divisor. The
  corresponding reflexive sheaf $\cO_X(K_X)$ is defined as
  $i_*(\Omega_U^{\dim X})$, where $i:U\to X$ is the embedding of the
  nonsingular part of $X$.
\end{say}

\begin{defn}
  An {\em $\Bbb R$-divisor\/} $D=\sum d_{j}D_{j}$ is a linear
  combination of prime Weil divisors with real coefficients, i.e.\ an
  element of $N^{1}\otimes \Bbb R$.  An $\Bbb R$-divisor is said to be
  $\Bbb R$-Cartier if it is a combination of Cartier divisors with
  real coefficients, i.e.\ if it belongs to the image of
  $Div(X)\otimes\Bbb R \to N^{1}(X)\otimes\Bbb R$ (this map is of
  course injective for normal varieties).  The $\Bbb Q$-divisors and
  $\Bbb Q$-Cartier divisors are defined in a similar fashion.
\end{defn}

\begin{defn}
  Consider an $\Bbb R$-divisor
  $K+B=K_{X}+\sum b_{j}B_{j}$ and assume that
  \begin{enumerate}
  \item $K+B$ is $\Bbb R$-Cartier,
  \item $0\le b_{j}\le1$.
  \end{enumerate}
  For any resolution $f:Y\to X$ look at the natural formula
  \begin{eqnarray}
    \label{defn:codiscrepancies}
    K_{Y}+B^{Y}= f^{*}(K_{X}+\sum b_{j}B_{j})=
    K_{Y}+\sum b_{j}f^{-1}B_{j} + \sum b_{i}F_{i}
  \end{eqnarray}
  or, equivalently,
  \begin{eqnarray}
    \label{defn:log discrepancies}
    K_{Y}+\sum b_{j}f^{-1}B_{j} + \sum F_{i}=
    f^{*}(K_{X}+\sum b_{j}B_{j}) + \sum a_{i}F_{i}
  \end{eqnarray}
    Here $f^{-1}B_{j}$ are the proper preimages of $B_{j}$ and $F_{i}$
  are the exceptional divisors of $f:Y\to X$.

  The coefficients $b_{i},b_{j}$ are called codiscrepancies, the
  coefficients $a_{i}=1-b_{i},a_{j}=1-b_{j}$ -- log discrepancies.
\end{defn}

\begin{rem}
  In fact, $K+B$ is not a usual $\Bbb R$-divisor but rather a special
  gadget consisting of a linear class of a Weil divisor $K$ (or a
  corresponding reflexive sheaf) and an honest $\Bbb R$-divisor $B$.
  This, however, does not cause any confusion.
\end{rem}

\begin{defn}
  A pair $(X,B)$ (or a divisor $K+B$) is said to be
  \begin{enumerate}
  \item log canonical, if the log discrepancies $f_{k}\ge0$
  \item Kawamata log terminal, if $f_{k}>0$
  \item canonical, if $f_{k}\ge1$
  \item terminal, if $f_{k}>1$
  \end{enumerate}
  for every resolution $f:Y\to X$, $\{k\}=\{i\} \cup \{j\}$.
\end{defn}

\begin{say}
  The notion of {\em semi-log canonical\/} is a generalization of
  {\em log canonical\/} to the case of varieties that are singular in
  codimension 1. The basic observation here is that for a curve with a
  simple node the definition of the log discrepancies still makes
  sense and gives $a_1=a_2=0$, so it can also be considered to
  be (semi-)log canonical. No new Kawamata semi-log terminal
  singularities appear, however.

  Recall that according to Serre's criterion normal is equivalent to
  Serre's condition $S_2$ and regularity in codimension 1. So, if
  we do allow singularities in codimension 1, $S_2$ will be exactly
  what we will need to keep.
\end{say}

\begin{defn}
  \label{defn:semi-log canonical}
  Let $X$ be a reduced (but not necessarily irreducible)
  equidimensional scheme which satisfies Serre's condition $S_{2}$ and
  is Gorenstein in codimension 1. Let $B=\sum b_{j}B_{j}$, $0\le
  b_{j}\le1$ be a linear combination with real coefficients of
  codimension 1 subvarieties none of irreducible components of which
  is contained in the singular locus of $X$.  Denote by $\cO(K_X)$ the
  reflexive sheaf $i_*(\omega_U)$, where $i:U\to X$ is the open subset
  of Gorenstein points of $X$ and $\omega_U$ is the dualizing sheaf of
  $U$. We can again consider a formal combination of $K_X$ and an
  $\Bbb R$-divisor $B$, and there is a good definition for $K_X+\sum
  b_jB_j$ to be $\Bbb R$-Cartier.  It means that in a neighborhood of
  any point $P\in X$ we can choose a section $s$ of $\cO(K_X)$ such
  that the divisor $(s)+\sum b_jB_j$ is a formal combination with real
  coefficients of Cartier divisors with no components entirely in the
  singular set.

  A pair $(X,B)$ (or a divisor $K_{X}+B$) is said to be semi-log
  canonical if, similar to the above,
  \begin{enumerate}
  \item $K_{X}+B$ is $\Bbb R$-Cartier,
  \item for any morphism $f:Y\to X$ which is birational on every
    irreducible component, and with a nonsingular $Y$, in the natural
    formula
    \begin{displaymath}
      f^*(K_{X}+B)=K_{Y}+f^{-1}B+ \sum b_{i}F_{i}
    \end{displaymath}
    with $F_{i}$ being irreducible components of the exceptional set,
    all $b_{i}\le 1$ (resp.~ $a_{i}=1-b_{i}\ge0$).
  \end{enumerate}
  As before, the coefficients $b_i,b_j$ are called codiscrepancies,
  the coefficients $a_i,a_j$ -- the log discrepancies.
\end{defn}

\begin{rem}
  In the case when $(X,B)$ has a good semi-resolution (for example,
  for surfaces) this definition is equivalent to that of
  \cite{KollarShepherdBarron88}, \cite{FAAT} chapter 12. In our
  opinion, it is more natural to give a definition which is
  independent of the existence of a semi-resolution.
\end{rem}

\begin{rem}
  For surfaces the condition $S_2$ is of course equivalent to
  Cohen-Macaulay.
\end{rem}

\begin{defn}
  By the Kleiman's criterion, the ampleness for proper schemes is a
  numerical condition, hence it extends to $\Bbb R$-Cartier divisors.
  If coefficients of $B$ are rational, $K_{X}+B$ is $g$-ample iff for
  some positive integer $n$ the divisor $n(K_X+B)$ is Cartier and
  $g$-ample in the usual sense.
\end{defn}

\begin{rem}
  Below we will only consider the case when $B$ is reduced, i.e.\ all
  the coefficients $b_j=1$. See the last section for the discussion on
  non-integral coefficients.
\end{rem}

\begin{exmp}
  If $X$ is a curve then $(X,B)$ is semi-log canonical iff the only
  singularities of $X$ are simple nodes and $B$ consists of distinct
  points lying in the nonsingular part of $X$.  $K_X+B$ is relatively
  ample iff every smooth rational component of $X$ mapping to a point
  on $W$ has at least 3 points of intersection with the rest of $X$, or
  points in $B_j$, and every component of arithmetical genus 1 has at
  least 1 such point. In the absolute case, i.e. when $W$ is a point,
  this is the usual definition of a stable curve with marked points.
  Every $B_j$ can also be considered as a group of unordered points.
\end{exmp}

\begin{exmp}
  The only codimension 1 semi-log canonical singularities are normal
  crossings.
\end{exmp}

\begin{exmp}
  If $X$ is a nonsingular surface then $(X,B)$ is semi-log canonical
  iff $B$ has only normal intersections.
\end{exmp}

\begin{exmp}
  For the case when $X$ is a surface and $B$ is empty the semi-log
  canonical singularities over $\Bbb C$ were classified in
  \cite{KollarShepherdBarron88}. They are (modulo analytic
  isomorphism): nonsingular points, Du Val singularities, cones over
  nonsingular elliptic curves, cusps or degenerate cusps (which are
  similar to cones over singular curves of arithmetical genus 1),
  double normal crossing points $xy=0$, pinch points $x^{2}=y^{2}z$,
  and all cyclic quotients of the above.  If $B$ is nonempty then the
  singularities of $X$ are from the same list and, in addition, there
  are different ways $B$ can pass through them. For normal $X$ the
  list could be found in \cite{Alexeev92} for example.
\end{exmp}

\begin{say}
  The following describes an easy reduction of semi-log canonical
  singularities to log canonical, cf. \cite{FAAT} 12.2.4.
\end{say}

\begin{lem}
  \label{lem:reduction of semi-log canonical to log canonical}
  Let $(X,B)$ be as in the definition \ref{defn:semi-log canonical}
  and denote by $\nu:X^{\nu}\to X$ its normalization.  Assume that
  $K_X+B$ is $\Bbb R$-Cartier.  Then $(X,B)$ is semi-log canonical iff
  $(X^{\nu},\nu^{-1}B+cond(\nu))$ is log canonical, and they have the
  same log discrepancies.
\end{lem}
\begin{pf}
  Clear from the definition.
\end{pf}

\begin{say}
  The next theorem explains how semi-log canonical surfaces appear in
  families (cf. \cite{KollarShepherdBarron88} 5.1).  But first we will
  need an auxiliary definition.
\end{say}

\begin{defn}
  Let $f:(\cX,\cB)\to \cS$ be a 3-dimensional one-parameter family.
  Let $\cB=\sum b_{j}\cB_{j}$ with $0\le b_{j} \le1$ be an
  $\Bbb R$-divisor and assume that $\cX$ and all $\cB_j$ are flat over
  $\cS$ and that $K_{\cX}+\cB$ is $\Bbb R$-Cartier.  We will say that
  the pair $(\cX,\cB)$ (or the divisor $K_{\cX}+\cB$) is {\em
    $f$-canonical} if in the definition of log discrepancies for all
  exceptional divisors with $f(F_i)$ a closed point on $\cS$ one has
  for the corresponding log discrepancy $a(F_i)\ge1$ (resp.
  $b(F_i)\le0$). This condition does not say anything about log
  discrepancies of divisors that map surjectively onto $\cS$.
\end{defn}

\begin{thm}
  \label{thm:family is good iff central fiber is good}
  Let $f:(\cX,\cB)\to \cS$ be a 3-dimensional one-parameter family
  over a pointed curve or a specter of a DVR (a discrete valuation
  ring).  Let $\cB=\sum b_{j}\cB_{j}$ with $0\le b_{j} \le1$ be an
  $\Bbb R$-divisor and assume that $\cX$ is irreducible, $\cX$ and all
  $\cB_j$ are flat over $\cS$ and that the fibers satisfy Serre's
  condition $S_2$ and are Gorenstein in codimension 1 (note that this
  implies that $\cX$ itself is Cohen-Macaulay and is Gorenstein in
  codimension 1).  Further assume that $K_{\cX}+\cB$ is $\Bbb
  R$-Cartier. Then the following is true:
  \begin{enumerate}
  \item If $K_{\cX_{0}}+\cB_{0}$ is semi-log canonical then
    $K_{\cX}+\cB$ is log canonical and $f$-canonical.
  \item Under assumptions of (1), the general fiber is also semi-log
    canonical.
  \item Suppose that there exists a birational morphism
    $\mu: \cY \to \cX$ with a nonsingular $\cY$ such that all
    exceptional divisors of $\mu$ and strict preimages of $\cB_{i}$
    have normal crossings and such that the central fiber is reduced.
    Then the opposite to (1) is true.
  \end{enumerate}
\end{thm}

\begin{pf}
  The proof is an application of the adjunction formula.

  (1) The log adjunction theorem \cite{FAAT} 17.12 and
  \ref{lem:reduction of semi-log canonical to log canonical} imply
  that $K_{\cX_0}+\cB_{0}$ is semi-log canonical iff
  $K_{\cX}+\cB+\cX_0$ is.  Now, the connection between the log
  discrepancies of the divisors $K_{\cX}+\cB$ and $K_{\cX}+\cB+\cX_0$
  is clear. For a divisor $E$ mapping onto $\cS$ the log discrepancies
  are the same.  For $E$ mapping to a central point of $\cS$ the
  difference is the coefficient of $E$ in the central fiber of
  $\cY\to\cS$, and so is at least 1.

  (3) Here the differences between the log discrepancies over the
  central fiber are all equal to 1.

  (2) is by adjunction.
\end{pf}

\begin{say}
  Finally, we show how to pass from a relatively ample $K+B$ to an
  ample divisor.
\end{say}

\begin{lem}
  \label{lem:absolute ampleness}
  Let $g:(X,B)\to W\subset \Bbb P$ be a map, where $X$ is a
  projective surface and $B=\sum b_jB_j$ is an $\Bbb R$-divisor on
  $X$.  Assume that $K_X+B$ is semi-log canonical and is relatively
  $g$-ample. Then $K_X+B+4H$ is ample, where $H=g^*\cO(1)$.
\end{lem}
\begin{pf}
  It is enough to prove that the restriction on the normalization of
  $X$ is ample, therefore by
  \ref{lem:reduction of semi-log canonical to log canonical} we can
  assume that $X$ is normal and that $(X,B)$ is log canonical.

  We show that $K_X+B+3H$ is nef (numerically effective) and this
  implies the statement. Indeed, $K_X+B+MH$ is ample for $M\gg0$ and
  $K_X+B+4H$ is a weighted average of the above two divisors.

  Assume that $K_X+B+3H$ is not nef. Then the Cone Theorem, which
  holds for arbitrary normal surfaces, tells us that there exists an
  irreducible curve $C$ generating an extremal ray and such that
  \begin{displaymath}
    (K_X+B)C<0
  \end{displaymath}
  This is possible only if $C$ does not map to a point. But then
  $C\cdot3H\ge3$ and $(K_X+B)C\ge-3$ by a theorem on the length
  of extremal curves, see \cite{MiyaokaMori86}.

  In dimension 2 the latter statement is very elementary.  Let
  $f:Y\to X$ be a minimal resolution of singularities of $X$. Then,
  if $X\ne\Bbb P^2$, one necessarily has $(f^{-1}C)^2\le0$ and
  \begin{eqnarray*}
    (K_X+B)C \ge (K_Y+f^{-1}B)f^{-1}C \ge \\
    (K_Y+f^{-1}C)f^{-1}C = 2p_a(f^{-1}C)-2 \ge -2.
  \end{eqnarray*}
  And the case of $X=\Bbb P^2$ is clear.
\end{pf}

\section{Definition and properties of a moduli functor}
\label{sec:definition and properties of a moduli functor}

\begin{say}
  Below we give a few general definitions for moduli functors. They
  are fairly standard (see e.g. \cite{Viehweg94}, \cite{Kollar90}) but
  we need to make slight modifications to adapt them to our situation.
\end{say}

\begin{say}
  The moduli functor for a moduli problem of polarized schemes is
  normally constructed in the following way. One fixes a class $\cC$
  of schemes $X/k$ with a polarization, i.e. an ample line bundle, $L$
  and some extra structure and subject to some conditions. Then for an
  arbitrary scheme $\cS/k$ one defines $\cMC(\cS)$ as the set of all
  (relatively) polarized flat families over $\cS$ with all geometric
  fibers from $\cC$ and, possibly, subject to more conditions. The
  families are considered modulo an equivalence relation. Usually it
  is an isomorphism between $\cX_1/\cS$ and $\cX_2/\cS$ with whatever
  extra structure they have and a fiber-wise linear equivalence
  between $\cL_1$ and $\cL_2$. In other cases it is an algebraic or a
  numerical, or a numerical up to a scalar equivalence, instead of
  linear.

  Sometimes, it is also useful considering a $\Bbb Q$-polarization
  $L$ on $X$, i.e a reflexive sheaf such that $(L^{\otimes N})^{**}$
  is an ample line bundle.
\end{say}

\begin{say}
  The above definition is intentionally vague since extra structures
  and conditions vary greatly from one moduli problem to another.
  Instead of trying to cover all future generalizations, we will
  formulate general principles and, when nontrivial, say exactly how
  they specialize to our situation.
\end{say}

\begin{defn}
  The class $\cS$ is said to be {\em bounded\/} if there exists a
  scheme $(\cX,\cL)$ with an extra structure, and a morphism
  $F:\cX\to\cS$ to a scheme $\cS$ of finite type such that all
  elements of $\cC$ appear as geometric fibers of $F$, not necessarily
  in a one-to-one way.

  There are two important variations of this definition. There is the
  {\em polarized\/} boundedness, when one requires $F$ to be
  projective and $\cL$ to restrict to the given polarization $L$ on a
  fiber, versus {\em non-polarized\/}. One can also consider
  boundedness {\em in the narrow sense\/}, requiring that all fibers
  of $F$ belong to $\cC$, or {\em in the wide sense\/}, asking only
  for some of the fibers to be from $\cC$.

  Here we make the choice of the polarized boundedness in the wide
  sense.
\end{defn}

\begin{defn}
  A moduli functor $\cMC$ is said to be {\em separated\/} if every
  one-parameter family in ${\cM\cC(\cS_{gen})}$, where $\cS_{gen}$ is
  a generic point of a DVR, has at most one extension to $\cS$.
\end{defn}

\begin{defn}
  A moduli functor $\cMC$ is said to be {\em complete\/} if every
  one-parameter family in ${\cM\cC(\cS_{gen})}$, where $\cS_{gen}$ is
  a generic point of a DVR, has at least one extension after a finite
  cover ${\cS'\to\cS}$.
\end{defn}

\begin{defn}
  \label{defn:local closedness}
  A class $\cC$ is said to be {\em locally closed\/} if for
  every flat family $F:(\cX,\cL)\to\cS$ with an extra structure
  there exist locally closed subschemes $\cS_l\subset\cS$ with the
  following universal property:
  \begin{itemize}
  \item A morphism of schemes $\cT\to\cS$ factors through
    $\coprod S_l$ iff $(\cX,\cL)\underset{\cS}{\times}\cT \to\cT$
    belongs to $\cMC(\cT)$.
  \end{itemize}
\end{defn}

\begin{defn}
  The class $\cC$ is said to {\em have finite reduced automorphisms\/}
  if every object in $\cC$ has a finite and reduced (the latter is
  automatic in characteristic 0) group of automorphisms.
\end{defn}

\begin{defn}
  A moduli functor $\cMC$ is said to be {\em functorially polarizable\/}
  if for every family $(\cX,\cL)$ in \cMCS there exists an equivalent
  family $(\cX,\cL^c)$ such that
  \begin{enumerate}
  \item if $(\cX_1,\cL_1)$ and $(\cX_2,\cL_2)$ are equivalent, then
    $(\cX_1,\cL^c_1)$ and $(\cX_2,\cL^c_2)$ are isomorphic,
  \item for any base chance $h:{\cS'\to\cS}$, $(\cX',{\cL'}^c)$ and
    $(\cX',h^*(\cL^c))$ are isomorphic.
  \end{enumerate}
  The main example of a functorial polarization is delivered by the
  polarization $\omega_{\cX/\cS}$ for canonically polarized manifolds.
\end{defn}

\begin{defn}
  A functorial polarization $\cL^c$ is said to be {\em semipositive\/}
  if there exists a fixed $k_{0}$ such that whenever $\cS$ is a
  complete smooth curve and $f:(\cX,\cL)\to\cS$ an element in \cMCS,
  then for all $k\ge k_{0}$ the vector bundles $f_*(k\cL^c)$ are
  semipositive, i.e.\ all their quotients have nonnegative degrees.

  This definition will be slightly modified for our purposes, we will
  also require semipositiveness of restrictions of $\cL^c$ to certain
  divisors $\cB_j$ on $\cX$.
\end{defn}

\begin{say}
  The following is the class that we will be considering from now on.
\end{say}

\begin{defn}
  \label{defn:the class}
  The elements of the class $\cC^N=\cC^N_{(K+B)^2,(K+B)H,H^2}$ are
  {\em stable maps of pairs\/} $g:(X,B,L_N)\to W$, where
  \begin{enumerate}
  \item $W\subset\Bbb P$ is a fixed projective scheme,
  \item $X$ is a connected projective surface,
  \item $B=\sum_{j=1}^n B_j$ is a divisor on $X$, $B_j$ are reduced
    but not necessarily irreducible,
  \item the pair $(X,B)$ has semi-log canonical singularities,
  \item the divisor $K_X+B$ is relatively $g$-ample,
  \item $(K_X+B)^2=C_1, (K_X+B)H=C_2, H^2=C_3$ are fixed,.
  \item $L_N=\cO(N(K_X+B+5H))$, where $H=g^*\cO_W(1)$. Here $N$ is a
    positive integer such that for every map as above $L_N$ is a line
    bundle. For example, we can choose $N$ to be the minimal positive
    integer satisfying this condition. The existence of such an $N$
    will be proved in \ref{thm:boundedness of maps}, and it is ample
    by \ref{lem:absolute ampleness}.
  \end{enumerate}
\end{defn}

\begin{say}
  The classes $\cC^N$ and $\cC^M$ for different $N,M$ are in a
  one-to-one correspondence between each other, and the only
  difference is the polarizations.

  As a consequence, the polarization in our functor plays a secondary
  role. We will switch from a polarization $L_N$ to its multiple $L_M$
  when it will be convenient.
\end{say}

\begin{thm}
  \label{thm:boundedness of maps}
  For some $M>0$ the class $\cC^M$ is bounded.
\end{thm}
\begin{pf}
  We start with the boundedness theorem which gives what we want in
  the absolute case.
  \begin{thm}[\cite{Alexeev94b}, 9.2]
    \label{thm:absolute boundedness}
    Fix a constant $C$ and a set $\cA$ satisfying the descending
    chain condition. Consider all surfaces $X$ with an $\Bbb R$-divisor
    $B=\sum b_jB_j$ such that the pair $(X,B)$ is semi-log canonical,
    $K_X+B$ is ample, $b_j\in \cA$ and $(K_X+B)^2=C$. Then the class
    $\{(X,\sum b_jB_j)\}$ is bounded.
  \end{thm}
  Apply this theorem with the set $\cA=\{1\}$ to $K_X+B+D$, where $D$
  is a general member of the linear system $|4H|$. Since this linear
  system is base point free, the pair $(X,B+D)$ also has semi-log
  canonical singularities.  Therefore, all pairs $(X,B)$ satisfying
  the conditions of the theorem can be embedded by a linear system
  $|M(K_X+B+4L)|$ for a fixed large divisible $M$ in a fixed
  projective space $\Bbb P^{d_1}$.  Every map $g:X\to W$ is defined by
  its graph $\Gamma_g$.  Consider a Veronese embedding of $W$ by
  $|\cO_W(M)|$ in some $\Bbb P^{d_2}$ and then look at the graphs
  $\Gamma_g$ in a Segre embedding
  $\Bbb P^{d_1}\times\Bbb P^{d_2} \subset \Bbb P^{d_3}$. Note that
  $\cO_{\Bbb P^{d_3}}(1)$ restricted on $X\simeq\Gamma_g$ is
  $L_M=M(K_X+B+5H)$.

  $L^2$ is fixed, hence by the boundedness theorem
  \ref{thm:boundedness of maps} above there are
  only finitely many possibilities for Hilbert polynomials
  $\chi(\cO_{\Gamma_g}(t))$. By the same theorem, there are also only
  finitely many possibilities for Hilbert polynomials
  $\chi(\cO_{B_j}(t))$. Therefore, all elements of our class
  $g:(X,B)\to W$ are parameterized by finitely many products of
  Hilbert schemes. In each product, we have to extract a subscheme
  parameterizing subschemes of $\Bbb P^{d_1}\times W$ and with fixed
  $\cO_{\Bbb P^{d_1}}(1)^2$,
  $\cO_{\Bbb P^{d_1}}(1)\cdot\cO_{\Bbb P^{d_2}}(1)$ and
  $\cO_{\Bbb P^{d_2}}(1)^2$, and these are obviously closed algebraic
  conditions. We also need to extract the graphs, i.e subschemes
  mapping isomorphically to $\Bbb P^{d_1}$, and this is an open
  condition.

  The resulting scheme will parameterize the maps, including all maps
  from the class $\cC^M$. This proves the theorem.
\end{pf}

\begin{say}
  We won't need the boundedness of the class $\cC^N$ itself, although
  it will follow from the proof of the local closedness
  \ref{thm:local closedness}.
\end{say}

\begin{defn}
  \label{defn:moduli functor 1}
  There are several ways to define the moduli functor for our class.
  The one we use here is the most straightforward one (cf.
  \cite{KollarShepherdBarron88}, \cite{Viehweg94} in the absolute case
  with $B=\emptyset$). For any scheme $\cS/k$,
  $\cMC^N=\cMC^N_{(K+B)^2,(K+B)H,H^2}$ is given by
  \begin{displaymath}
    {\cMC^N(\cS)}=
    \left\{
      \begin{aligned}
        & \text{all families }
        f:(\cX,\cL)\to\cS
        \text{ with a divisor } \cB=\sum_{j=1}^N \cB_j
        \text{  on } X, \\
        & \text{a map }
        g:\cX\to W
        \text{ and a line bundle } \cL
        \text{ such that every }\\
        & \text{geometric fiber belongs to } \cC,
        X \text{ and all } \cB_j
        \text{ are flat over } \cS \\
      \end{aligned}
    \right\}
  \end{displaymath}
  Two families over $\cS$ are equivalent if they are isomorphic
  fiber-wise.

  In this functor we consider a sub-functor $\cMC{'}^N$, requiring
  in addition that for each $s$ there exists a
  1-dimensional family from $\cMC^N$ with the central fiber
  $\cX_s$ and an irreducible general fiber $\cX_g$ such that:
  \begin{enumerate}
  \item $\cX_g$ is irreducible,
  \item the pair $(\cX_g,0)$ is (Kawamata) log terminal.
  \end{enumerate}

  This is similar to the smoothability condition for
  $\overline{M_{K^2}^{sm}}\subset\overline{M_{K^2}}$ (see
  \cite{Kollar90}) and is necessary due to the technical reasons.
  Consider a one parameter family of maps. Then we would like the
  ambient 3-fold to be irreducible since MMP is not developed for
  non-irreducible varieties yet. We would also want the 3-fold to have
  log terminal singularities because they are Cohen-Macaulay in
  characteristic 0.
\end{defn}

\begin{say}
  A little disadvantage of the above definition is that even though
  $\cMC^{N,irr}$ and, say, $\cMC^{2N,irr}$ are the same on the closed
  points, the corresponding moduli spaces can potentially have
  different scheme structures, the second one could be strictly
  larger.  So, in fact, we have not one but infinitely many moduli
  spaces.  It would be better if we had a formula for the minimal $N$
  in terms of $(K+B)^2,(K+B)H,H^2$. We know, however, only that such
  an $N$ exists.
\end{say}

\begin{say}
  A different solution was suggested (again, in the absolute case with
  $B=\emptyset$) by Koll\'ar in \cite{Kollar90},\cite{Kollar94}. In a
  sense, it produces a moduli space with the ``minimal'' scheme
  structure.

  We introduce some necessary notation first.
\end{say}

\begin{defn}
  Let $F:\cX\to\cS$ be a projective family of graphs of maps
  $(X,B)\to W$. Assume that every fiber is Gorenstein in codimension 1
  and satisfies Serre's condition $S_2$.
  Denote by $i:\cU\hookrightarrow\cX$ the open subset where $f$ is
  Gorenstein and the divisors $\cB_j$ are Cartier.  Note that on every
  fiber one has $\codim_{\cX_s}(\cX_s-\cU_s)\ge2$. Define the sheaves
  $\cL_{\cU,k}$ and $\cL_k$ by
  \begin{displaymath}
    \cL_{\cU,k}=\cO_{\cU}(k(K_{\cU/\cS}+\cB+g^*\cO_W(5))
  \end{displaymath}
  \begin{displaymath}
    \cL_k=i_*\cL_{\cU,k}
  \end{displaymath}
  It follows that the sheaves $\cL_k$ on $\cX$ are coherent.
\end{defn}

\begin{notationnum}
  Let $f:\cX\to\cS$ be a morphism of schemes,
  $i:\cU\hookrightarrow\cX$ be the immersion of an open set and $\cF$ be
  a coherent sheaf on $\cU$ which is flat over $\cS$. For a base
  change $h:\cS'\to\cS$ we obtain
  $\cX^h:=\cX\underset{\cS}{\times}\cS'$,
  $\cU^h:=\cU\underset{\cS}{\times}\cS'$ etc. Denote the induced
  morphism $\cU^h\to \cU$ by $h_{\cU}$ and set $\cF^h:=h_{\cU}^*\cF$.
  The induced morphism $\cX^h\to\cX$ is denoted by $h_X$

  One says that
  {\em the push forward of $\cF$ commutes with a base change\/}
  $h:\cS'\to\cS$ if the natural map $h_X^*(i_*\cF)\to i^h_*\cF^h$ is an
  isomorphism.
\end{notationnum}

\begin{defn}
  Define $\cMC^{all}=\cMC^{all}_{(K+B)^2,(K+B)H,H^2}$ by
  \begin{displaymath}
    {\cMC^{all}(\cS)}=
    \left\{
      \begin{aligned}
        & \text{all families }
        f:\cX\to\cS
        \text{ with a divisor } \cB=\sum_{j=1}^N \cB_j
        \text{ on } \cX
        \text{ and} \\
        & \text{a map }
        g:\cX\to W
        \text{ such that every geometric fiber belongs to }
        \cC, \\
        & \cX \text{ and all } \cB_j
        \text{ are flat over } \cS,
        \text{ and for each } k \\
        & i_*\cL_{\cU,k}
        \text{ commutes with arbitrary base changes}
      \end{aligned}
    \right\}
  \end{displaymath}
  As above, one can consider a sub-functor
  $\cMC{'}^{all}\subset\cMC^{all}$.

  We will not go into detailed discussion of this functor.
\end{defn}

\begin{say}
  One can see that if we require that $i_*\cL_{\cU,k}$ commutes with
  arbitrary base changes only for $k=N$ instead of all positive $k$,
  then we get the previous definition of the moduli functor.  Indeed,
  if a line bundle $\cL$ exists, then $\cL_N=\cL+f^*\cE$ for some
  invertible sheaf $\cE$ on $\cS$. Then for every $h:\cS'\to\cS$ the
  two sheaves $i^h_*\cL_{\cU,N}^h$ and
  $h_X^*(i_*\cL_{\cU,N})=h_X^*(\cL_N)$ on $\cX'$ are both reflexive
  and coincide on $h_X^{-1}(\cU)$, hence everywhere.

  Vice versa, if $i_*\cL_{\cU,N}$ commutes with base changes, then
  $\cL_N$ is flat and for every closed point $s\in\cS$
  \begin{displaymath}
    \cL_N\big|_{\cX_s}=\cO_{\cX_s}(N(K+B+g^*\cO_W(5)))
  \end{displaymath}
  Since the latter restriction is locally free for every $s$ and the
  sheaves $\cO_{\cX}$, $\cL_N$ are coherent and flat over $\cS$, it
  follows by \cite{Matsumura86} 22.5, 22.3 that $\cL_N$ is locally
  free.
\end{say}

\begin{say}
  Now let us show that our moduli functor $\cMC{'}^N$ has all the good
  properties listed above. We start with the local closedness. The
  main technical result we will be using is the following theorem.
\end{say}

\begin{thm}[Koll\'ar \cite{Kollar94}]
  \label{kollar's flattening decomposition}
  With the above notations, assume that $f:\cX\to\cS$ is projective,
  $i_*\cF$ is coherent and that for every point $s\in\cS$ the sheaf
  $\cF_s$ on the fiber $\cX_s$ satisfies Serre's condition $S_2$.
  Then there exist locally closed subschemes $\cS_l\subset\cS$ such
  that for any morphism $h:\cT\to\cS$ the following are equivalent:
  \begin{enumerate}
  \item $h$ factors through $\cT\to\coprod\cS_l\to\cS$,
  \item $i^h_*\cF^h$ commutes with all future base changes.
  \end{enumerate}
\end{thm}

\begin{thm}
  \label{thm:local closedness}
  The functors $\cMC^N$ and $\cMC{'}^N$ are locally closed.
\end{thm}
\begin{pf}
  Let $F:\cX\to\cS$ be an arbitrary projective family of graphs of
  maps $(X,B)\to W$. First, after the flattening decomposition (see
  \cite{Mumford66} lecture 8) of $\cS$ into locally closed subschemes,
  we can assume that $\cX$ and $\cB_j$ are flat over $\cS$ if they are
  not already.

  Consider a one-parameter sub-family $\cX_{\cR}\to\cR$ and a point $P$
  on the central fiber $\cX_0$. Then $\cX_0$ is Cohen-Macaulay at $P$
  iff the 3-fold $\cX_{\cR}$ is. The property of a
  local ring to be Cohen-Macaulay
  is open (\cite{Matsumura86} 24.5) and the morphism $F$ is
  projective. Therefore, if $\cX_0$ is Cohen-Macaulay then there exists
  an open neighborhood of $\cR$, and also of $\cS$, that contains
  exactly the points over which the fibers are Cohen-Macaulay.

  The property of a local ring to be Gorenstein is also open
  (\cite{Matsumura86} 24.6) and by the same argument there exists a
  closed subset $Z$ of non-Gorenstein points in $\cX$.  Give it the
  structure of a reduced scheme. Then we have to throw away all fibers
  on which the Hilbert polynomial of
  $\cO_Z\underset{\cO_{\cS}}{\otimes}k(s)$ has degree $\ge1$. There
  are only finitely many possible Hilbert polynomials and the
  condition on the degree is obviously closed.

  At this point we use the previous theorem
  \ref{kollar's flattening decomposition} to the sheaf $\cL_{\cU,N}$
  to conclude that there exist locally closed subschemes
  $\cS_l\subset\cS$ such that every map $h:\cT\to\cS$ with
  $\cX\underset{\cS}{\times}{\cT}\in\cMC(\cT)$ factors through
  $\coprod\cS_l$. $\cS_l$ are disjoint, so
  we can concentrate on one of them. If $P$ is a point of $\cS$
  and some $h$ as in the definition does not factor through $\cS-P$,
  then the fiber of $F$ over $P$ has to be a pair  $(X,B)$ from our
  class. The sheaf $\cL_N$ on $X\underset{\cS}{\times}S_l$ is flat
  over $S_l$ and its restriction to the fiber over $P$ is locally
  free. Hence, it has to be locally free in a neighborhood of the
  fiber. Therefore, for each $S_l$ if we denote by $U_l\subset S_l$
  the open set over which $\cL_N$ is locally free, then $h:\cT\to\cS$
  has to factor through $\coprod\cU_l$. Now we can apply
  \ref{thm:family is good iff central fiber is good}(2) to conclude that
  there exist open subsets $\cV_l\subset\cU_l$ containing all
  the points over which the fibers have semi-log canonical
  singularities.  Also, $\cMC{'}^N\subset\cMC^N$ is evidently closed
  and we end up with a disjoint union of locally closed subschemes.

  There is one more thing one has to take care of: the polarization
  $\cO_{\Bbb P^{d_3}}(1)$ on the fibers has to coincide with $\cL_N$
  or its fixed multiple $\cL_M$. Standard semi-continuity theorems for
  $h^0$ in flat families show that there exists a closed subset where
  the two sheaves are the same. One can also define the scheme
  structure on it, see lemma 1.26 \cite{Viehweg94}.
\end{pf}

\begin{lem}
  For the functors $\cMC^N$ and $\cMC{'}^N$ the polarization $\cL_N$ is
  functorial.
\end{lem}
\begin{pf}
  $K_{\cU/\cS}$ of a flat family commutes with base changes, and so
  do $\cO(\cB_j)$ and $g^*\cO_W(1)$. Therefore, $\cL_{\cU,k}$ are
  functorial.  By the definitions of the functor $\cMC$ the same is
  true for $\cL_k$ (resp. $\cL_N$).
\end{pf}

\begin{thm}
  $\cMC{'}^N$ is
  \begin{enumerate}
  \item separated,
  \item complete,
  \item have finite and reduced automorphisms.
  \end{enumerate}
\end{thm}
\begin{pf}
  The first two properties have code names in the Minimal Model
  Program: ``uniqueness and existence of the log canonical model''.
  It is enough to check them in the case when the general fiber is
  irreducible and has log terminal singularities.

  (1) Let $\cS$ be a specter of a DVR or a pointed curve. Two families
  in $\cMC(\cS)$ that coincide outside of $0$ are birationally
  isomorphic.  \ref{thm:family is good iff central fiber is good}(1)
  implies that they are both log canonical and both are relative log
  canonical models over $\cS\times W$ for the same divisor, hence
  isomorphic. If $\cY\to\cS$ is a common resolution then the divisor
  is
  \begin{displaymath}
    K_{\cY}+f^{-1}\cB + \sum\cE_i
  \end{displaymath}
  where $\cE_i$ are exceptional divisors that do not map to a central
  point $0\in\cS$.

  (2) If there is a family over $\cS-0$, we can complete it over $0$
  somehow. Then by a variant of the Semistable Reduction Theorem,
  after a finite base change, there is a resolution $\cY$ of
  singularities such that the central fiber is reduced and all
  exceptional divisors and $\cB_j$ have normal crossings. Consider the
  log canonical model for the same divisor as above, relative over
  $\cS\times W$. It exists by \cite{KeelMcKernanMatsuki93} for
  example.  This log canonical model has the same fibers as $(X,B)$
  outside $0$. It has log terminal
  singularities only, which are Cohen-Macaulay in dimension 3 and
  characteristic 0. Therefore, the central fiber is also Cohen-Macaulay
  and it is from our class $\cC$ by
  \ref{thm:family is good iff central fiber is good}(3).

  We also have to show that the sheaf $\cL_N$ for this family is
  locally free. It amounts to proving that the Hilbert polynomials
  $h_1(t)$ of the sheaf $L_{N,0}$ on the special fiber, and $h_2(t)$
  of the sheaf $L_{N,g}$ of the general fiber coincide. Both sheaves
  are locally free.  But the log canonical model is constructed by
  applying the Base Point Freeness theorem, and by the very
  construction we have that some $\cL_M$ for a large divisible $M$
  is locally free on $\cX$. Therefore the polynomials $h_1(M/Nt)$ and
  $h_2(M/Nt)$ are the same,
  and that means that $h_1(t)$ and $h_2(t)$ are also the same.

  (3) In the absolute case, the fact that $K+B$ is ample and log
  canonical implies that the automorphism group is finite by
  \cite{Iitaka82}. In the relative case we apply the same theorem to
  $K_X+B+D$, $D\in|4H|$ general, which is ample by lemma
  \ref{lem:absolute ampleness}.  We are working in characteristic 0
  and so the group scheme $\Aut X$ is reduced.
\end{pf}

\begin{thm}
  \label{thm:semipositiveness}
  The functors $\cMC^N$ and $\cMC{'}^N$ are semipositive.
\end{thm}
\begin{pf}
  One has the following
  \begin{thm}[Koll\'ar \cite{Kollar90} 4.12]
    \label{thm:kollar's semipositiveness}
    Let $Z$ be a complete variety over a field of characteristic zero.
    Assume that $Z$ satisfies Serre's condition $S_2$ and that it is
    Gorenstein in codimension one. Let $Z\to C$ be a map onto a smooth
    curve. Assume that the general fiber of $f$ has only semi-log
    canonical singularities, and further that $K$ of the general fiber
    is ample. Then $f_*\cO(kK_{Z/C})$ is semipositive for $k\ge1$.
  \end{thm}

  For the sheaves $\cL_N= O_{\cX}(N(K_{\cX/\cS}+{\cB}))$ with empty
  $\cB$ in the absolute case this is exactly what we need. Analyzing
  the proof of \ref{thm:kollar's semipositiveness} shows that it works
  with very minor changes in the case of a non-empty reduced $\cB$. In
  the relative case instead of $K_{\cX/\cS}+{\cB}$ we consider
  $K_{\cX/\cS}+{\cB}+5H$, $H=g^*\cO_W(1)$. We can think of $5H$ simply
  as of an additional component of the boundary $\cB$.  If a member of
  the linear system $|5H|$ is chosen generically, on the general fiber
  of $f$ the pair $(X,B+5H)$ will still be semi-log canonical.

  For the positiveness of the sheaves $\cL_N\Big|_{\cB_j}$ we use the
  log adjunction formula, see \cite{Shokurov91} or \cite{FAAT} chapter
  16. We get the following semi-log canonical divisors on $\cB_j$:
  \begin{displaymath}
    K_{\cX}+\cB\Big|_{\cB_j}=K_{\cB_j}+\sum(1-1/m_{k})\cM_{k}
  \end{displaymath}
  for some Weil divisors $\cM_{k}$ on $\cB_j$ and $m_k\in\Bbb
  N\cup\{\infty\}$. So, here we need a more general semipositiveness
  theorem, with nonempty $\cB$ that has fractional coefficients. The
  situation is saved by the fact that the relative dimension of
  $\cB_j$ over $\cS$ equals 1, and the semipositiveness for this case
  is proved in \cite{Kollar90} 4.7.
\end{pf}

\section{Existence and projectivity of a moduli space}
\label{sec:existence and projectivity of a moduli space}

\begin{thm}
  \label{thm:existence as an algebraic space}
  The functor $\cMC=\cMC{'}^N$ is coarsely represented by a proper
  separated algebraic space of finite type $\bMC=\bMC{'}^N$.
\end{thm}
\begin{pf}
  The proof is essentially the same as in \cite{Mumford82}, p.172. We
  remind that we are working in characteristic zero, and over $\Bbb C$
  the argument is easier.

  The class $\cC^M$ is bounded, and we can embed all graphs $\Gamma_g$
  of the maps $g$ by a linear system $|M(K_X+B+5H)|$ in $\Bbb
  P^{d_1}\times\Bbb P^{d_2} \subset \Bbb P^{d_3}$ as in
  \ref{thm:boundedness of maps} for a large divisible $M$.  By taking
  $M$ even larger we can assume that all $X=\Gamma_g$ and all
  $B_j\subset X$ are projectively normal, $h^0(M(K_X+B+5H))$ is
  locally constant and there are no higher cohomologies.
  $(\Gamma_g,B)$ are parameterized, not in a one-to-one way, by some
  scheme that we will denote by $\cH$. For any family in
  $\cMC^N(\cT)$, the embedding by a relatively very ample linear
  system $|M(K_X+B+4H)|$ defines a non-unique map $\cT\to\cH$.  By
  \ref{thm:local closedness} there exists a disjoint union of locally
  closed subschemes $\cS=\coprod\cS_l\hookrightarrow\cH$ with a
  universal property, and $\cT\to\cH$ factors through $\cS$.  We
  conclude that the coarse moduli space $\bMC$ is a categorial
  quotient of $\cS$ by an equivalence relation $R$, described as
  follows.

  $R$ is a set of pairs $(h,G)$, where $h\in\cS$ and $G$ corresponds
  to a different embedding of $X$ in $\Bbb P^{d_1}$, i.e.  $G$ varies
  in a group $PGL(d_1+1)$. There is a natural map
  $F:R\to\cS\times\cS$.
  Every fiber of $\pi_1\circ F$ is isomorphic to $PGL(d_1+1)$ and this
  map is obviously smooth. The map $F$ is quasi-finite and unramified
  because its fibers are automorphism groups of objects in $\cC$,
  and these are finite reduced. The fact that $\cMC$ is also proper
  implies that $F$ is finite.

  The rest of the proof is the same as in \cite{Mumford82}, p.172
  verbatim. By taking the transversal sections locally the question is
  reduced to the case of a finite equivalence relation dominated by a
  map $F':R\to\cH'\times\cH'$ with $\pi_{1}\circ F'$ \'etale, and then the
  quotient is easily constructed as an algebraic space.

  Finally, since $\cMC$ is proper, so is $\bMC$.
\end{pf}

\begin{thm}
  \label{thm:projectiveness}
  The moduli space $\bMC=\bMC{'}^N$ is projective.
\end{thm}
\begin{pf}
  The proof follows the general scheme of \cite{Kollar90}.  By the
  very construction of $\bMC$, there exists a subscheme
  $\cS\subset\cH$ of a product of Hilbert schemes, with the
  corresponding universal family $V_{\cS}\to\cS$, that maps to $\bMC$.
  One starts by constructing a {\em finite\/} morphism from a scheme
  $Y\to\bMC$ with a universal family $f:V_Y\to Y$.  This is done
  locally by cutting $\cS\to\bMC$ transversally, then adding more
  copies of these sections, so that the automorphisms do not obstruct
  gluing the local pieces together, see \cite{Kollar90} 2.7.  The only
  properties of the class $\cC$ used in this construction are
  boundedness and finiteness of automorphisms, which we have.

  Next step is to consider the line bundles
  \begin{displaymath}
    \lambda_M=\det(f_*\cL_M\oplus f_*\cL_M\big|_{\cB_j})
  \end{displaymath}
  on $Y$ for $M$ large divisible, where
  \begin{displaymath}
    \cL_M=\cO_V(M(K_{V/Y}+\cB+g^*\cO_W(5))).
  \end{displaymath}

  These line bundles do not descend to $\bMC$ because of
  automorphisms, but since the objects of $\cC$ have finite groups of
  automorphisms and $\cC$ is bounded, for every $M$ there is a finite
  power of $\lambda_M$ that does come from a line bundle on $\bMC$. To
  prove that $\bMC$ is projective it is enough to show that one of
  $\lambda_M$ is ample, which is achieved by the following theorem.
  For simplicity we formulate it only in characteristic 0.

  \begin{notationnum}
    Let $Y$ be a scheme and let $W$ be a vector bundle of rank $w$
    with structure group $\rho:G\to GL_w$. Let $q:W\to Q$ be a
    quotient vector bundle of rank $k$. Let $Gr(w,k)/G$ denote the set
    of $G$-orbits on the $k$-dimensional quotients of a
    $w$-dimensional vector space. The natural map of sets
    \begin{displaymath}
      u_{Gr}:\{\text{closed points of }X\}\to
      Gr(w,k)/G
    \end{displaymath}
    is called the {\em classifying map}.

One says that the classifying map is {\em finite\/} if
\begin{enumerate}
\item every fiber of $u_{Gr}$ is finite, and
\item for every $y\in Y$ only finitely many elements of $G$ leave
  $\ker q_y$ invariant.
\end{enumerate}
  \end{notationnum}
  \begin{thm}[Koll\'ar's Ampleness Lemma, \cite{Kollar90} 3.9]
    Let $Y$ be a proper algebraic space and let $W$ be a semipositive
    vector bundle with structure group $G$. Let $Q$ be a quotient
    vector bundle of $W$. Assume that
    \begin{enumerate}
    \item $G$ is reductive,
    \item the classifying map is finite.
    \end{enumerate}
    Then $\det Q$ is ample. In particular, $Y$ is projective.
  \end{thm}
  This is what it translates to in our situation. The sheaves are
  \begin{displaymath}
    W=\Sym^j(f_*\cL_M)\oplus\Sym^j(f_*\cL_M\big|_{\cB_j})
  \end{displaymath}
  and
  \begin{displaymath}
    Q=f_*L_{jM}\oplus f_*L_{jM}\big|_{\cB_j},
  \end{displaymath}
  $q$ is the multiplication map.

  By \ref{thm:semipositiveness} we already know that $Q$ is
  semipositive, and so is $W$ since symmetric powers of a semipositive
  sheaf are semipositive.

  Recall that the universal family $U_Y$ over $Y$ is embedded into a
  product of $Y$ and
  \begin{displaymath}
    \Bbb P^{d_1}\times W \subset
    \Bbb P^{d_1}\times\Bbb P^{d_2} \subset \Bbb P^{d_3}
  \end{displaymath}
  and that the sheaf $L_{M}$ is the restriction of
  $\cO_{P^{d_3}}(1)$ in this embedding.

  The group $G$ acting on $W$ is $GL_{d_1+1}\times GL_1$.  If every
  fiber $\cX=\cGamma_g$ together with all $B_j$ can be uniquely
  reconstructed from the map $W_s\to Q_s$, then the fibers of $u_{Gr}$
  will be exactly the same as fibers of $Y\to\bMC$, hence finite.

  For this to be true we need the following:
  \begin{enumerate}
  \item every fiber in $\Bbb P^{d_3}$ is set-theoretically defined by
    degree $\le j$ equations,
  \item the multiplication maps $\Sym^j(f_*\cL_M)\to f_*L_{jM}$
    and $\Sym^j(f_*\cL_M\big|_{\cB_J})\to f_*L_{jM}\big|_{\cB_J}$
    are surjective.
  \end{enumerate}
  (1) holds if $j$ is large enough. (2) is satisfied because we have
  chosen $M$ so large that all $X$ and $B_j$ are projectively
  normal in $\Bbb P^{d_3}$.

  Finally, the second condition in the definition of finiteness of the
  classifying map is satisfied because all graphs $(\Gamma_g,B)=(X,B)$
  in $\Bbb P^{d_3}$ have finite groups of automorphisms.
\end{pf}

\section{Related questions}
\label{sec:related questions}

\begin{say}
  Let us see how our moduli spaces are related to some others.  For
  example, consider the moduli space $\cM_{L^2}$ of K3 surfaces $X$
  with a polarization $L$ with a fixed square.
  Compare it with $\cM_{(K+B)^2}$, where $W=pt$, $B=B_1$ is one
  reduced divisor and $(K+B)^2=L^2$ is the same number.
  $\cM_{(K+B)^2}$ contains an open subset $U$ parameterizing K3
  surfaces with reduced divisors having normal intersections only, and
  we have a map $F:U\to\cM_{H^2}$. A well-known result (Saint-Donat
  \cite{SaintDonat74}) says that every ample linear system $|L|$ on a
  K3 surface contains at least one reduced divisor with normal
  intersections, therefore $F$ is surjective. In fact, $\cM_{H^2}$ is
  a quotient of $U$ modulo an obvious equivalence relation $R$:
  $(X_1,B_1)\underset{R}{\sim}(X_2,B_2)$ iff $X_1$, $X_2$ are
  isomorphic and $B_1$, $B_2$ are linearly equivalent.

  There is a natural map $G:R\to U\times U$.  $\pi_1\circ G$ is smooth
  and its fibers are open subsets in $\Bbb P^{h^0(H)-1}$. The
  situation is very similar to what we had in theorem
  \ref{thm:existence as an algebraic space}, except this time the
  quotient $U/R$ is not proper. The obvious way to try to obtain a
  compactification of $\cM_{H^2}$ is to consider the closure
  $\overline{U}$ of $U$ in $\cM_{(K+B)^2}$, then somehow define
  the closure $\overline{R}$ of $R$, and ask if it has good enough
  properties enabling one to construct $\overline{U}/\overline{R}$ and
  to prove that it is projective. Alternatively, one can ask if the
  closure of $\overline{G}(R)$ in $\overline{U}\times\overline{U}$ has
  good properties.

  The situation resembles what happens for elliptic curves. The
  natural compactification of the moduli space $\cM_1=\Bbb A^1_k$ is
  $\Bbb P_k^1$, and the infinite point corresponds not to one but to
  many degenerations: wheels of rational curves of lengths $1\dots n$
  if we consider $\cM_1$ as a factor of $\cM_{1,n}$. Similarly, the
  boundary points of $\cM_{H^2}$ should correspond to many different
  degenerations of smooth K3 surfaces with geometric divisors,
  properly identified.

  The first thing to ask on this way is:
  \begin{question}
    Is it possible to define an equivalence relation
    $\overline{G}:\overline{R}\to\overline{U}\times\overline{U}$,
    so that the morphism $\pi_1\circ\overline{G}$ is smooth or at
    least flat?
  \end{question}

  Even if this is done, there are problems with taking the quotient.
  There does not seem to exist in the literature a ready-to-use method
  that would cover our situation. There is, on one hand, a theorem of
  M.Artin (see \cite{Artin69} 7.1, \cite{Artin74b} 6.3) that shows that
  if $\overline{G}:\overline{R}\to\overline{U}\times\overline{U}$ were
  a  monomorphism (which it is not) with flat projections, then the
  quotient would be defined as an algebraic space. In this case it
  would also easily follow that the quotient is actually projective.

  On the other hand, there is the method of \cite{Mumford82}, p.172
  that we used in the previous section, in which the equivalence
  relation is smooth, and the map $\overline{G}$ is finite. Natural
  degenerations of K3 surfaces can have infinite groups of
  automorphisms, however.

  I think that the question deserves a more detailed consideration.
\end{say}

\begin{say}
  Similarly to K3 surfaces, for any principally polarized Abelian
  variety $A$ with a theta divisor $\Theta$ the pair $(A,\Theta)$ has
  log canonical singularities, see \cite{Kollar93}. So. the previous
  discussion applies to principally polarized Abelian surfaces too.
  One can also ask what happens if the polarization is not principal.
\end{say}

\begin{say}
  It goes without saying that the projectivity theorem
  \ref{thm:projectiveness} applies in the case of curves, with
  significant simplifications.  Therefore, the moduli spaces $M_{g,n}(W)$
  of \cite{Kontsevich94} are also projective.
\end{say}

\begin{say}
  Most $\bMC_{(K+B)^2,(K+B)H,H^2}$ are definitely not irreducible and
  not even connected. They are subdivided according to various
  invariants, such as the numerical or homological type of $g(X)$ and
  $g(B_j)$, intersection numbers $(K+B)B_j$ etc.

  One can also get by fixing only one number, $(K+B+4H)^2$. Then there
  are only finitely many possibilities for other invariants.
\end{say}

\begin{say}
  The boundedness theorem \ref{thm:absolute boundedness} is in fact
  even stronger than what we used here: it applies to the case when
  the coefficients $b_j$ belong to an arbitrary set $\cal A$ that
  satisfies the descending chain condition. One, perhaps, would want
  to define even more general moduli spaces.  There are two obstacles,
  however.  First, the semipositiveness theorem \ref{thm:kollar's
    semipositiveness} for the case of fractional coefficients seems to
  be quite hard to prove, but probably still possible. The second
  obstacle is a fundamental one: for proving the semipositiveness
  theorems for $\cL_k|_{\cB_j}$ we used the log adjunction formula. It
  basically just says $K+B|_B=K_B$, and here the coefficient 1 of $B$
  is important.
\end{say}

\begin{say}
  The places where assumption about the characteristic 0 was
  used:
  \begin{enumerate}
  \item MMP in dimension 3. This is not serious since we worked in the
    situation of the relative dimension 2. For surfaces log MMP is
    characteristic free, and perhaps it is true for families of
    surfaces in generality needed. For the case $B=\emptyset$ see
    \cite{Kawamata91}
  \item The semipositiveness theorem
    \ref{thm:kollar's semipositiveness} requires characteristic 0.
    Since we are dealing with a case of relative dimension 2 only,
    this also probably can be dealt with.
  \item A group scheme in characteristic 0 is reduced, hence
    smooth. This was used in the proof of
    \ref{thm:existence as an algebraic space}. Perhaps, the argument
    could be strengthened.
  \item The argument of \cite{Mumford82} p.172 is a whole lot more
    complicated in characteristic $p>0$.
  \end{enumerate}
\end{say}

\begin{say}
  It should be possible to prove the semipositiveness theorems and the
  Ampleness Lemma, as well as the \cite{Mumford82} p.172 argument,
  entirely in the relative situation $/W$, without appealing to
  absolutely ample divisors. The moduli spaces obtained should be then
  projective over $W$.
\end{say}

\begin{say}
  One can see that most of the theorems that we proved for the functor
  $\cMC{'}^N$ apply to the functor $\cMC{'}^{all}$ as well.
\end{say}

\makeatletter \renewcommand{\@biblabel}[1]{\hfill#1.}\makeatother

\end{document}